\documentclass[12pt]{article}  
\usepackage{amsmath,amssymb}

\usepackage{graphicx}

\usepackage{color}

\newdimen\mypaperwidth
\newdimen\mypaperheight
\def\mya4{\mypaperwidth=210mm
\mypaperheight=297mm}
\def\mycenteralloc{%
\oddsidemargin=\mypaperwidth
\advance\oddsidemargin-\textwidth
\oddsidemargin=0.5\oddsidemargin
\advance\oddsidemargin-1in
\evensidemargin=\oddsidemargin
}
\mycenteralloc
\def\dsp{\displaystyle}

\def\U#1{\underline{#1}}

\def\1{{\mathbf 1}}

\def\figtxt#1#2#3{\ignorespaces\rlap{\kern #1\raisebox{#2}{#3}}\ignorespaces}

\def\surd{\sqrt{\vrule width0pt height1.2ex}}

\def\bar#1{\setbox0\hbox{$#1$}\dimen0=\ht0\advance\dimen0by0.8pt
  \overline{\vrule width0pt height\the\dimen0 #1\kern0.25pt}}
\def\db#1{\widetilde{#1\kern1pt}}

\def\bar#1{\overline{#1}}

\author{Shigekazu Nakagawa, Hiroki Hashiguchi$^{\dag}$ and Yoko Ono$^{\ddag}$}
\date{Okayama University of Science, Tokyo University of Science$^{\dag}$, and Yokohama City University$^{\ddag}$}

\title{Approximation to probability density functions in sampling distributions based on Fourier cosine series}

\begin{document}
\maketitle

\begin{abstract}
We derive a simple and precise approximation to probability density functions in sampling distributions based on the Fourier cosine series.
After clarifying the required conditions, we illustrate the approximation on two examples: the distribution of the sum of uniformly distributed random variables, and   
the distribution of sample skewness drawn from a normal population.
The probability density function of the first example can be explicitly expressed, but that of the second example has no explicit expression.
\end{abstract}

\section{Introduction}

Fourier series and their transforms are widely applied in statistics because they are mathematically tractable. For instance, they are uniformly convergent within a closed interval and can be integrated
term by term (Whittaker~\cite{whittaker1963course}).
Substantial applications are 
the generalized representations (or estimations) of statistical curves such as
probability density functions (pdfs) and their associated cumulative distribution functions (cdfs).
(see, for example, Kronmal and Tarter~\cite{KronmalTarter1968}).

Many approximations to pdfs or cdfs have been developed in sampling distribution theory.
The Edgeworth series of the pdf of a statistic, a well-known series that
refines the central limit theorem, is composed of Hermite polynomials as orthogonal polynomials.
Meanwhile, the Fourier series are composed of orthogonal trigonometric functions and their approximation target slightly differs from that of the Edgeworth series.
Whereas the Edgeworth series handles cases with bounded of unbounded pdf support, 
the Fourier series requires a bounded support and a periodic pdf.
It appears that the Fourier coefficients must be obtained under severe constraints.

In this paper, we approximate the pdf or cdf using cosine Fourier series, and
 provide three conditions that must be satisfied.
First, the pdf requires a bounded support.
Second, the function must be piecewise smooth and even.
Finally, the pdf must have moments up to the required order. 
We remark that the first and second conditions are technical only, but the third is crucial.
We also remark that sampling distributions are evaluated by numerous statistics, for example, 
sample skewness, sample kurtosis, the Shapiro--Wilk test statistic, and sample correlation coefficient.

The remainder of this paper is organized as follows. 
Section~\ref{sec:bs} provides a brief survey of Fourier cosine series, and 
Sections~\ref{sec:uniform} and \ref{sec:skewness} demonstrate our approach on two examples.
More specifically, Section~\ref{sec:uniform} considers the distribution of the sum of uniformly distributed random variables, for which the pdf can be explicitly expressed, and Section~\ref{sec:skewness} considers a statistic with no explicit pdf expression but with moments 
up to the required order. 
The accuracy of the proposed method is evaluated in numerical experiments of the first case. 
In the second case, the distribution of sample skewness $\surd{b_{1}}$ is drawn from a normal population. For samples of size $n = 3$ and $4$, 
the pdfs were obtained by Fisher~\cite{Fisher1930} and McKay~\cite{mckay1933distribution}, respectively; for other sample sizes, they are unknown. 
Here, we mention the sample of size $n$. 
When $n$ is large, $\surd{b_{1}}$ is asymptotically normally distributed with zero mean and variance $n/6$ (Thode~\cite{thode2002testing}).
When $n$ is moderate,  D'Agostino's~\cite{d1970transformation} transformation works well.
When $n$ is at most 25, Mulholland's~\cite{mulholland1977null} approximations are valid, 
but the mathematical expressions of the pdfs are very complicated.

Based on the Fourier cosine series, 
we provide concrete approximations of the sampling distribution, which are particularly effective when $n$ is small.
We also give the percentiles of $\surd{b_{1}}$.

\section{Fourier series of probability density functions}
\label{sec:bs}

Let $\dsp T_{n} = T_{n}\left( X_{1}, X_{2}, \ldots, X_{n} \right)$ be a 
statistic with probability density function $f_{n}(x)$, where 
$\left( \dsp X_{1}, X_{2}, \ldots, X_{n} \right)$ is a random sample of size $n$.
We assume that the following conditions are satisfied:
\begin{itemize}
\item[1.]
$f_{n}(x)$ has a bounded support $[-A_{n}, A_{n}]$, where $A_{n} > 0$.
\item[2.]
$f_{n}(x)$ is both piecewise smooth and even.
\end{itemize}

\noindent
From Condition 1 and 2, 
$f_{n}(x)$ has a Fourier cosine series within $[-A_{n}, A_{n}]$.
That is,
\begin{equation}
\label{eq:f_pdf}
f_{n}(x) = \frac{a_{n, 0}}{2} + 
\sum_{k=1}^{\infty} 
a_{n, k} \cos \frac{k\pi}{A_{n}}x, 
\end{equation}
where 
the Fourier cosine coefficients $a_{n, k}$ are given by
\begin{align}
a_{n, 0} & = \frac{1}{A_{n}} \int_{-A_{n}}^{A_{n}} f_{n}(x) \, dx = \frac{1}{A_{n}},
\nonumber
\\
a_{n, k} & = \frac{1}{A_{n}} \int_{-A_{n}}^{A_{n}} f_{n}(x) \cos \frac{k\pi}{A_{n}}x \, dx.
\label{eq:ak_org}
\end{align}

\noindent
We also assume that the moments of $T_{n}$ about the origin 
exist up to the requisite order:
\begin{itemize}
\item[3.]
For any $j$,  there exists $\mu_{n, j}^{\prime} = \dsp \int_{-\infty}^{\infty} x^j f_{n}(x) \, dx  < \infty$.
\end{itemize}

\noindent
From Condition 3 and 
$\dsp\cos x = \sum_{j=0}^{\infty}\frac{(-1)^{j}}{(2j)!} x^{2j}$, 
the coefficient 
$a_{n, k}$ is determined as 
\begin{align}
a_{n, k}
& =
\frac{1}{A_{n}} \int_{-A_{n}}^{A_{n}} f_{n}(x) \sum_{j=0}^{\infty} \frac{(-1)^{j}}{(2j)!} \left( \frac{k\pi}{A_{n}} \right)^{2j} x^{2j} \, dx
\nonumber
\\
& =
\frac{1}{A_{n}} 
\sum_{j=0}^{\infty} \frac{(-1)^{j}}{(2j)!} \left( \frac{k\pi}{A_{n}} \right)^{2j} \mu_{n, 2j}^{\prime},
\label{eq:ak_org_2}
\end{align}
where we have performed termwise integration.

The cumulative distribution function $F_{n}(x) = \Pr\left\{ T_{n} < x \right\}$ 
is also obtained as
\begin{equation*}
\label{eq:f_cdf}
F_{n}(x) = \frac{1}{2}\left( \frac{x}{A_{n}} + 1 \right) 
+ \sum_{k=1}^{\infty}  \frac{a_{n, k}A_{n}}{k\pi} \sin \frac{k\pi}{A_{n}} x.
\end{equation*}

\medskip
From the Fourier cosine series \eqref{eq:f_pdf} of $f_{n}(x)$, we find that the approximations
\begin{equation}
\widetilde{f_{n}}^{(K)}(x) = 
\frac{a_{n, 0}}{2} + 
\sum_{k=1}^{K} 
a_{n, k} \cos \frac{k\pi}{A_{n}}x
\label{eq:pdf_approx}
\end{equation}
and
\begin{equation}
\widetilde{F_{n}}^{(K)}(x) = \frac{1}{2}\left( \frac{x}{A_{n}} + 1 \right) 
 + \sum_{k=1}^{K}\frac{a_{n, k}A_{n}}{k \pi} \sin\frac{k\pi}{A_{n}} x 
\label{eq:cdf_approx}
\end{equation}
are the best approximations because they minimize
$$\dsp \int_{-A_{n}}^{A_{n}} \left( \widetilde{f_{n}}^{(K)}(x) - f_{n}(x) \right)^2\, dx \qquad \mbox{and} \qquad 
\dsp \int_{-A_{n}}^{A_{n}} \left( \widetilde{F_{n}}^{(K)}(x) - F_{n}(x) \right)^2\, dx,$$
respectively (Whittaker~\cite{whittaker1963course}).
\begin{itemize}
\item
Solving the nonlinear equation $\widetilde{F_{n}}^{(K)}(x) = \alpha$
by Newton's method with an initial value of $x_{0.5}=0$, the percentile is obtained as $x_{\alpha} \, (0 < \alpha < 1)$.
\item
If $f_{n}(x)$ is given, the coefficients $a_{n,k}$ are directly obtained from \eqref{eq:ak_org}.
If $f_{n}(x)$ is not given but its moments $\mu_{n, j}^{\prime}$ are known to the requisite order, 
the approximation equation \eqref{eq:pdf_approx} can be obtained from \eqref{eq:ak_org_2}.
%
\end{itemize}

\smallskip
We approximate the Fourier cosine coefficients \eqref{eq:ak_org} as
\begin{equation*}
\widehat{a_{n, k}}^{(J)} = 
\frac{1}{A_{n}} 
\sum_{j=0}^{J} \frac{(-1)^{j}}{(2j)!} \left( \frac{k\pi}{A_{n}} \right)^{2j} \mu_{n, 2j}^{\prime}
\end{equation*}
and Eqs. \eqref{eq:pdf_approx} and \eqref{eq:cdf_approx} as
\begin{align}
\widehat{f_{n}}^{(K,J)}(x) & = \frac{a_{n, 0}}{2} + 
\sum_{k=1}^{K}
\widehat{a_{n, k}}^{(J)} \cos \frac{k\pi}{A_{n}} x, 
\nonumber
\\
\widehat{F_{n}}^{(K,J)}(x) &= \frac{1}{2}\left( \frac{x}{A_{n}} + 1 \right) 
 + \sum_{k=1}^{K}\frac{\widehat{a_{n, k}}^{(J)}A_{n}}{k \pi} \sin\frac{k\pi}{A_{n}} x, 
\nonumber
\end{align}
respectively.
Note that the choices of $K$ and $J$ depend on $n$.

\section{Distribution of a sum of random variabels uniformly distributed}
\label{sec:uniform}

Random variables $X_{1}, X_{2},  \ldots , X_{n}$ are mutually independent 
and uniformly distributed on the interval $\dsp \left[ -\frac{1}{2}, \, \frac{1}{2}\right]$.
Let $f_{n}(x)$ be the probability density function of 
the sum $T_{n} = X_{1} + X_{2} +  \cdots + X_{n}$.
Using the relation 
$$
f_{n}(x)
 = \int_{-\infty}^{\infty} f_{n-1}(x-t)f_{1}(t) \, dt
 = \int_{x-1/2}^{x+1/2} f_{n-1}(t) \, dt
 \qquad (n \geq 2),
$$
$$
f_{1}(x) = 
\begin{cases}
    1  & \left( -\frac{1}{2} \leq x \leq \frac{1}{2} \right) \\
    0      &    (\rm{otherwise}), 
\end{cases}
$$
for any 
$\dsp -\frac{n}{2} + i \leq x \leq  -\frac{n}{2} + i +1 \, (i= 0, 1, \ldots , n -1)$, we have
$$
f_{n}(x) = \frac{1}{(n-1)!} \sum_{j=0}^{i}(-1)^{j}{n \choose j}\left(x + \frac{n}{2} - j \right)^{n-1}, 
$$
$f_{n}(x) = 0 \, (\mbox{otherwise})$.
That is, $f_{n}(x)$ is both smooth and even, and has the bounded support $\dsp \left[ -\frac{n}{2}, \frac{n}{2} \right]$.
The moments $\mu_{n, j}^{\prime} $ are given by the following recurrence formula:
$$
\mu_{n, 2j}^{\prime}
= \sum_{k=0}^{j}{2j \choose 2k}\frac{1}{(2k+1)4^{k}}
\mu_{n-1, 2j-2k}^{\prime}.
$$

\bigskip
{\bf Example}
\quad When $n=4$, the probability density function is 
$$
f_{4}(x) = 
\begin{cases}
 \frac{1}{6} \left( x+2 \right)^3 
 & \left( -2 \leq x \leq -1  \right) \\
 \frac{1}{6} \left\{ \left(x+2 \right)^3 -4 \left(x+ 1 \right)^3 \right\}
 & \left( -1 \leq x \leq 0 \right) \\
 \frac{1}{6} \left\{  \left(x + 2\right)^3 -4 \left(x+1 \right)^3 +6 x^3 \right\}
 & \left( 0 \leq x \leq 1  \right) \\
 \frac{1}{6} \left\{  \left(x + 2\right)^3 -4 \left(x+1 \right)^3 +6 x^3  -4 \left(x-1 \right)^3 \right\}
 & \left( 1 \leq x \leq 2  \right) \\
    0      &    (\mbox{otherwise}).
\end{cases}
$$
The even moments are given by
$$\dsp \mu_{4, 2j}^{\prime} = \frac{8 \left( 4\cdot 4^j -1 \right)}
{(1 + 2j)(2+2j)(3+2j)(4+2j)}, $$ and 
odd moments are all $0$.

The Fourier cosine series of $f_{4}(x)$ in $[-2,2]$ is
\begin{equation*}
f_{4}(x) = \frac{1}{4} + 
\sum_{k=1}^{\infty}
a_{4, k} \cos \frac{k\pi}{2} x, 
\label{eq:f4_unif}
\end{equation*}
with coefficients
\begin{equation*}
a_{4, k} = \frac{128}{\pi^4  k^4} \sin ^4 \left( \frac{k\pi}{4} \right).
\label{eq:ak4_unif}
\end{equation*}

\noindent
The cumulative distribution function is
\begin{equation*}
F_{4}(x) = \frac{1}{2} + \frac{x}{4} 
+ \sum_{k=1}^{\infty}\frac{2a_{n, k}}{k\pi} \sin\frac{k\pi}{2} x. 
\label{eq:F4_unif}
\end{equation*}

\bigskip

Table~\ref{tab:ak} lists the
$$
\dsp \max_{0 \leq k \leq K} \left| a_{n, k}^{(J)} - a_{n, k}\right|
$$
values for $n = 2 (2) 12$. In all cases, the approximate coefficient $a_{n, k}^{(J)}$ coincided with the exact coefficient $a_{n, k}$ 
to six decimal places.

Tables~\ref{tab:ak48} and \ref{tab:ak812} tabulate the approximate Fourier cosine coefficients.
For example, 
\begin{align*}
\widehat{f_{4}}^{(8, 35)}(x) \, = & \,\,
   2.5\times 10^{-1}
+\left(3.28511\times 10^{-1}\right) \cos \left(\frac{\pi  x}{2}\right)+\left(8.21279\times
   10^{-2}\right) \cos (\pi  x)
   \\ &
   +\left(4.0557\times 10^{-3}\right) \cos \left(\frac{3 \pi 
   x}{2}\right)+\left(1.0306\times 10^{-14}\right) \cos (2 \pi  x)
   \\ &
   +\left(5.25618\times
   10^{-4}\right) \cos \left(\frac{5 \pi  x}{2}\right)
   +\left(1.01392\times 10^{-3}\right) \cos
   (3 \pi  x)
   \\ &
   +\left(1.36823\times 10^{-4}\right) \cos \left(\frac{7 \pi 
   x}{2}\right)-\left(5.73436\times 10^{-10}\right) \cos (4 \pi  x)
\end{align*}
and 
\begin{align*}
\widehat{F_{4}}^{(8, 35)}(x) \, = & \,\,
\frac{1}{2} \left(\frac{u}{2}+1\right)+\left(2.09137\times 10^{-1}\right) \sin \left(\frac{\pi 
   x}{2}\right)+\left(2.61421\times 10^{-2}\right) \sin (\pi  x)
   \\ &
   +\left(8.60646\times
   10^{-4}\right) \sin \left(\frac{3 \pi  x}{2}\right)+\left(1.64026\times 10^{-15}\right) \sin
   (2 \pi  x)
   \\ &
   +\left(6.69238\times 10^{-5}\right) \sin \left(\frac{5 \pi 
   x}{2}\right)+\left(1.07581\times 10^{-4}\right) \sin (3 \pi  x)
   \\ &
   +\left(1.24434\times
   10^{-5}\right) \sin \left(\frac{7 \pi  x}{2}\right)-\left(4.56326\times 10^{-11}\right) \sin
   (4 \pi  x).
\end{align*}

Table~\ref{tab:u4tail} lists the percentiles $x_{\alpha}$ of the sum $T_{n}$ for specified $\alpha$ values.
For example, when $n = 4$ and $\alpha = 0.99$, solving
$$
\widehat{F_{4}}^{(8, 35)}(x) = 0.99
$$
gives $x_{0.99} = 1.3002$.
We thus confirm that 
$$
\int_{-2}^{x_{0.99}} f_{4}(x) \, dx =  0.990006.
$$

\begin{table}[htbp]
\caption{Values of $\dsp \max_{0 \leq k \leq K} \left| a_{n, k}^{(J)} - a_{n, k}\right|$ for various $n$, $K$, and $J$}
\label{tab:ak}
\begin{center}
\begin{tabular}{rrrrr}
$n$ & $K$ & $J$ & $\dsp \max_{0 \leq k \leq K} \left| a_{n, k}^{(J)} - a_{n, k}\right|$
\\
\hline
2 & 8 & 35 & $3.61470 \times 10^{-7}$
\\
4 & 8 & 35 & $5.73436 \times 10^{-10}$
\\
6 & 8 & 30 & $1.13062 \times 10^{-7}$
\\
8 & 8 & 30 & $1.58854 \times 10^{-9}$
\\
10 & 7 & 25 & $6.93824 \times 10^{-10}$
\\
12 & 6 & 20 & $1.02801 \times 10^{-9}$
\end{tabular}
\end{center}
\end{table}

\begin{table}[htbp]
\caption{Fourier coefficients $\widehat{a}_{2, k}^{(35)},\, \widehat{a}_{4, k}^{(35)},\, \widehat{a}_{6, k}^{(30)}\, \, (k = 0, 1, \ldots , 8)$}
\label{tab:ak48}
$$
\begin{array}{rrrr}
k & \widehat{a}_{2, k}^{(35)}  &  \widehat{a}_{4, k}^{(35)} & \widehat{a}_{6, k}^{(30)}
\\ 
\hline
0& 1. & 5.\times 10^{-1} & 3.33333\times 10^{-1} \\
1& 4.05285\times 10^{-1} & 3.28511\times 10^{-1} & 2.52759\times 10^{-1} \\
2& 2.18614\times 10^{-16} & 8.21279\times 10^{-2} & 1.06633\times 10^{-1} \\
3& 4.50316\times 10^{-2} & 4.0557\times 10^{-3} & 2.21901\times 10^{-2} \\
4& 1.64157\times 10^{-13} & 1.0306\times 10^{-14} & 1.66614\times 10^{-3} \\
5& 1.62114\times 10^{-2} & 5.25618\times 10^{-4} & 1.61766\times 10^{-5} \\
6& 1.83356\times 10^{-11} & 1.01392\times 10^{-3} & -7.88258\times 10^{-15} \\
7& 8.27112\times 10^{-3} & 1.36823\times 10^{-4} & 2.14845\times 10^{-6} \\
8& -3.6147\times 10^{-7} & -5.73436\times 10^{-10} & 2.61465\times 10^{-5} \\
\end{array}
$$
\end{table}

\begin{table}[htbp]
\caption{Fourier coefficients $\widehat{a}_{8, k}^{(30)},\, \widehat{a}_{10, k}^{(25)},\, \widehat{a}_{12, k}^{(20)}\, \, (k = 0, 1, \ldots , 8)$}
\label{tab:ak812}
$$
\begin{array}{rrrr}
k & \widehat{a}_{8, k}^{(30)}  &  \widehat{a}_{10, k}^{(25)} & \widehat{a}_{12, k}^{(20)}
\\ 
\hline
0&  2.5\times 10^{-1} & 2.\times 10^{-1} & 1.66667\times 10^{-1} \\
1&  2.03319\times 10^{-1} & 1.69572\times 10^{-1} & 1.45271\times 10^{-1} \\
2&  1.0792\times 10^{-1} & 1.02664\times 10^{-1} & 9.58308\times 10^{-2} \\
3&  3.57613\times 10^{-2} & 4.34405\times 10^{-2} & 4.72705\times 10^{-2} \\
4&  6.74499\times 10^{-3} & 1.23309\times 10^{-2} & 1.70558\times 10^{-2} \\
5&  6.00653\times 10^{-4} & 2.18691\times 10^{-3} & 4.34421\times 10^{-3} \\
6&  1.64487\times 10^{-5} & 2.13837\times 10^{-4} & 7.38603\times 10^{-4} \\
7&  3.52694\times 10^{-8} & 9.08017\times 10^{-6} & \mbox{--} \\
8&  1.58854\times 10^{-9} & \mbox{--} & \mbox{--} \\
\end{array}
$$
\end{table}

\begin{table}[htbp]
\caption{Percentiles $x_{\alpha}$ of the sum $T_{n}$ for various $n$ and $\alpha$}
\label{tab:u4tail}
$$
\begin{array}{rrrrrrr}
n & \multicolumn{6}{c}{\alpha}
\\
\hline
  & 0.900 & 0.950 & 0.975 & 0.990 & 0.995 & 0.999 
\\
\hline
2& 0.5528 & 0.6838 & 0.7768 & 0.8571 & 0.8993 & 0.9649 \\
4& 0.7534 & 0.9534 & 1.1198 & 1.3002 & 1.4114 & 1.6063 \\
6& 0.9170 & 1.1663 & 1.3759 & 1.6097 & 1.7618 & 2.0536 \\
8& 1.0556 & 1.3457 & 1.5916 & 1.8694 & 2.0527 & 2.4120 \\
10& 1.1781 & 1.5039 & 1.7815 & 2.0971 & 2.3067 & 2.7232 \\
12& 1.2889 & 1.6469 & 1.9532 & 2.3028 & 2.5355 & 2.9964 \\
\end{array}
$$
\end{table}

\section{Distribution of a sample skewness drawn from normal population}
\label{sec:skewness}

Let $\left( X_{1}, X_{2}, \ldots, X_{n} \right)$ be a random sample of size $n$ drawn from a normal population.
Define
the sample skewness $\surd{b_{1}}$ as
$$
\surd{b_{1}} = \frac{m_{3}}{m_2^{3/2}}, \qquad
m_{r} = \frac{1}{n}\sum_{i=1}^{n} \left( X_{i} - \bar{X} \right)^{r}\,\, (r=2, 3), \qquad 
\bar{X} = \frac{1}{n}\sum_{i=1}^{n}X_{i}
$$
and let $f_{n}(x)$ be the probability density function.

Dalen~\cite{Dalen87} derived the range of $\surd{b_{1}}$ as
$$\displaystyle -A_{n} \leq \surd{b_{1}}  \leq A_{n}, \qquad A_{n} = \frac{n-2}{\sqrt{n-1}}.$$ 
Obviously, $f_{n}(x)$ is smooth and even on $[ - A_n , A_n ]$.
Geary~\cite{geary1947frequency} gave 
the recurrence formula of $f_{n}(x)$ as 
\begin{align}
f_{n}(x) & = \cfrac{\left(\frac{n-1}{n}\right)^{1/2}}{B\left(\frac{1}{2}, \frac{n-2}{2} \right)}
\int_{-1}^{1}
f_{n-1} \left( \sigma_{n-1}\left( x,z \right) \right) \left(1-z^2 \right)^{(n-7)/2} \, dz,
\label{eq:fnrec}
\end{align}
where 
$$\dsp \sigma_{n-1}\left( x, z \right)  = 
\left\{
\sqrt{n-1} x -3z + (n+1)z^3
\right\} n^{-1/2}\left( 1-z^2 \right)^{-3/2}.$$
However, an analytical expression of $f_{n}(x)$ is still lacking.
The moments of $\surd{b_{1}}$ are obtainable by the well-known 
recurrence relation obtained by Muholland~\cite{mulholland1977null}:
\begin{align}
\mu_{n+1, 2s}^{\prime} =& \cfrac{(n+1)^s}{n^s\left( \frac{n}{2}\right)_{3s}}
\sum_{j=0}^{s}{2s \choose 2j} \frac{\mu^{\prime}_{n, 2s-2j}}{(n+1)^j}
\label{eq:momentsrec}
\\
& \times
\sum_{i=0}^{2j}{2j \choose i} 3^{2j-i} (1-n)^i \left( \frac{1}{2} \right)_{j+i} 
\left( \frac{n-1}{2} \right)_{3s-j-i},
\nonumber
\end{align}
where $(a)_{m}$ is the Pochhammer symbol defined by 
$$
(a)_{m} = a(a+ 1) \cdots (a+m -1) \quad (m \geq 1), \qquad 
(a)_{0} = 1.
$$

\begin{figure}[htbp]
\begin{center}
\includegraphics{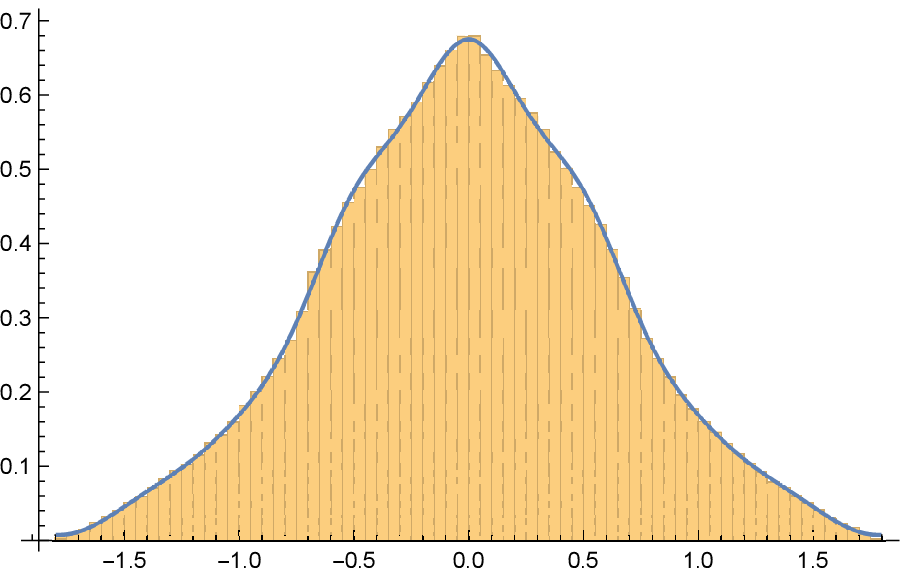}
\end{center}
\caption{Plot of $y=\widehat{f_{6}}^{(12, 50)}(x)$ and 
a histgram of $\surd{b_{1}}$ ($10^6$ replications)}
\label{fig:n6.pdf}
\end{figure}

From \eqref{eq:f_pdf}, the Fourier cosine series is obtained as 
\begin{equation}
\label{eq:sb1_pdf}
f_{n}(x) = \frac{a_{n, 0}}{2} + 
\sum_{k=1}^{\infty} a_{n, k} \cos \frac{k\pi}{A_{n}}x 
\end{equation}
and 
\begin{align*}
a_{n, 0}  = \frac{1}{A_{n}},
\quad
a_{n, k}  = \frac{1}{A_{n}}  \sum_{j=0}^{\infty} \frac{(-1)^j}{(2j)!} \left( \frac{k \pi}{A_{n}} \right)^{2j} \mu_{n,2j}^{\prime}.
\label{eq:ak}
\end{align*}
Using termwise integration, the cumulative distribution 
function $F_{n}(x) = \Pr\left\{  \surd{b_{1}} < x \right\}$ is obtained as 
\begin{equation*}
\label{eq:sb1_cdf}
F_{n}(x) = \frac{1}{2}\left( \frac{x}{A_{n}} + 1 \right) 
+ \sum_{k=1}^{\infty} a_{n, k} \frac{A_{n}}{k\pi} \sin \frac{k\pi}{A_{n}} x.
\end{equation*}
\noindent
Tables~\ref{tab:upptail410}, \ref{tab:upptail1218} and \ref{tab:upptail2022} show the 
approximate Fourier cosine coefficients $\widehat{a_{n, k}}^{(12, 50)} \, (k =0, \ldots, 12)$ for $n = 4 (2) 22$.
Especially, $f_{6}(x)$ and $F_{6}(x)$ are respectively approximated as 
\noindent
\begin{align*}
\widehat{f_{6}}^{(12, 50)}(x) \, = & \,\,
2.79508\times 10^{-1}
   + \left(3.08052\times 10^{-1}\right) \cos \left(\frac{1}{4} \sqrt{5} \pi 
   x\right)
   \\  &
   +\left(5.75070\times 10^{-2}\right) \cos \left(\frac{1}{2} \sqrt{5} \pi 
   x\right)+\left(1.17190\times 10^{-2}\right) \cos \left(\frac{3}{4} \sqrt{5} \pi 
   x\right)
   \\  &
   -\left(5.99392\times 10^{-3}\right) \cos \left(\sqrt{5} \pi 
   x\right)+\left(6.78429\times 10^{-3}\right) \cos \left(\frac{5}{4} \sqrt{5} \pi 
   x\right)
   \\  &
   +\left(7.71527\times 10^{-3}\right) \cos \left(\frac{3}{2} \sqrt{5} \pi 
   x\right)+\left(6.95419\times 10^{-3}\right) \cos \left(\frac{7}{4} \sqrt{5} \pi 
   x\right)
   \\  &
   +\left(1.62249\times 10^{-4}\right) \cos \left(2 \sqrt{5} \pi 
   x\right)+\left(2.38820\times 10^{-5}\right) \cos \left(\frac{9}{4} \sqrt{5} \pi 
   x\right)
   \\  &
   +\left(6.33581\times 10^{-4}\right) \cos \left(\frac{5}{2} \sqrt{5} \pi 
   x\right)+\left(3.10573\times 10^{-3}\right) \cos \left(\frac{11}{4} \sqrt{5} \pi 
   x\right)
   \\  &
   +\left(1.7351\times 10^{-3}\right) \cos \left(3 \sqrt{5} \pi 
   x\right)
   \end{align*}
and 
\begin{align*}
\widehat{F_{6}}^{(12, 50)}(x) \, = & \,\,
\frac{1}{2} \left(\frac{\sqrt{5} x}{4}+1\right)+\left(1.75408\times 10^{-1}\right) \sin
   \left(\frac{1}{4} \sqrt{5} \pi  x\right)
 \\ &
   +\left(1.63725\times 10^{-2}\right) \sin
   \left(\frac{1}{2} \sqrt{5} \pi  x\right)+\left(2.22431\times 10^{-3}\right) \sin
   \left(\frac{3}{4} \sqrt{5} \pi  x\right)
 \\ &
-\left(8.53250\times 10^{-4}\right) \sin
   \left(\sqrt{5} \pi  x\right)+\left(7.72609\times 10^{-4}\right) \sin \left(\frac{5}{4}
   \sqrt{5} \pi  x\right)
 \\ &
      +\left(7.32192\times 10^{-4}\right) \sin \left(\frac{3}{2}
   \sqrt{5} \pi  x\right)+\left(5.65684\times 10^{-4}\right) \sin \left(\frac{7}{4}
   \sqrt{5} \pi  x\right)
 \\ &   
   +\left(1.15483\times 10^{-5}\right) \sin \left(2 \sqrt{5} \pi 
   x\right)+\left(1.51096\times 10^{-6}\right) \sin \left(\frac{9}{4} \sqrt{5} \pi 
   x\right)
 \\ &   
   +\left(3.60767\times 10^{-5}\right) \sin \left(\frac{5}{2} \sqrt{5} \pi 
   x\right)+\left(1.60767\times 10^{-4}\right) \sin \left(\frac{11}{4} \sqrt{5} \pi 
   x\right)
    \\ &
   +\left(8.233\times 10^{-5}\right) \sin \left(3 \sqrt{5} \pi  x\right).
\end{align*}
Figure~\ref{fig:n6.pdf} superimposes the plot $y=\widehat{f_{6}}^{(12, 50)}(x)$ with 
a histogram of $\surd{b_{1}}$ (obtained after $10^6$ replications).

\begin{table}[htbp]
\caption{Fourier coefficients $\widehat{a_{n, k}}^{(50)}\quad \quad (n=4, 6, 8, 10)$}
\label{tab:upptail410}
\begin{center}
$$
\begin{array}{rrrrr}
k & n=4 & n = 6 & n=8 & n= 10 \\
\hline
0& 8.66025\times 10^{-1} & 5.59017\times 10^{-1} & 4.40959\times 10^{-1} & 3.75000\times
   10^{-1} \\
1& 1.76257\times 10^{-1} & 3.08052\times 10^{-1} & 3.12106\times 10^{-1} & 2.97971\times
   10^{-1} \\
2& 7.26283\times 10^{-2} & 5.75070\times 10^{-2} & 1.18300\times 10^{-1} & 1.55596\times
   10^{-1} \\
3& 5.56174\times 10^{-2} & 1.17190\times 10^{-2} & 3.20582\times 10^{-2} & 6.04639\times
   10^{-2} \\
4& 3.80806\times 10^{-2} & -5.99392\times 10^{-3} & 6.29481\times 10^{-3} & 1.95825\times
   10^{-2} \\
5& 3.28048\times 10^{-2} & 6.78429\times 10^{-3} & 1.54291\times 10^{-3} & 5.48500\times
   10^{-3} \\
6& 2.57761\times 10^{-2} & 7.71527\times 10^{-3} & 4.84745\times 10^{-4} & 1.49687\times
   10^{-3} \\
7& 2.32437\times 10^{-2} & 6.95419\times 10^{-3} & -5.28008\times 10^{-4} & 2.86985\times
   10^{-4} \\
8& 1.94759\times 10^{-2} & 1.62249\times 10^{-4} & -2.89091\times 10^{-4} & 5.48249\times
   10^{-5} \\
9& 1.79941\times 10^{-2} & 2.38820\times 10^{-5} & 6.47841\times 10^{-4} & 7.54897\times
   10^{-5} \\
10& 1.56490\times 10^{-2} & 6.33581\times 10^{-4} & 1.01858\times 10^{-3} & 2.55678\times
   10^{-5} \\
11& 1.46777\times 10^{-2} & 3.10573\times 10^{-3} & 6.95872\times 10^{-4} & -7.05180\times
   10^{-5} \\
12& 1.33900\times 10^{-2} & 1.73510\times 10^{-3} & 3.76670\times 10^{-4} & -5.51200\times
   10^{-5} \\
\end{array}
$$
\end{center}
\end{table}

\begin{table}[htbp]
\caption{Fourier coefficients $\widehat{a_{n, k}}^{(50)}\quad (n=12, 14, 16, 18)$}
\label{tab:upptail1218}
\begin{center}
$$
\begin{array}{rrrrr}
k & n=12 & n = 14 & n=16 & n= 18 \\
\hline
0& 3.31662\times 10^{-1} & 3.00463\times 10^{-1} & 2.76642\times 10^{-1} & 2.57694\times
   10^{-1} \\
1& 2.81214\times 10^{-1} & 2.65297\times 10^{-1} & 2.50979\times 10^{-1} & 2.38295\times
   10^{-1} \\
2& 1.75649\times 10^{-1} & 1.85442\times 10^{-1} & 1.89288\times 10^{-1} & 1.89688\times
   10^{-1} \\
3& 8.64878\times 10^{-2} & 1.06921\times 10^{-1} & 1.21861\times 10^{-1} & 1.32313\times
   10^{-1} \\
4& 3.61873\times 10^{-2} & 5.33942\times 10^{-2} & 6.92058\times 10^{-2} & 8.27429\times
   10^{-2} \\
5& 1.34095\times 10^{-2} & 2.39311\times 10^{-2} & 3.56438\times 10^{-2} & 4.73779\times
   10^{-2} \\
6& 4.55324\times 10^{-3} & 9.85304\times 10^{-3} & 1.69758\times 10^{-2} & 2.52441\times
   10^{-2} \\
7& 1.43298\times 10^{-3} & 3.79153\times 10^{-3} & 7.58376\times 10^{-3} & 1.26687\times
   10^{-2} \\
8& 4.17297\times 10^{-4} & 1.37582\times 10^{-3} & 3.21053\times 10^{-3} & 6.04348\times
   10^{-3} \\
9& 1.25392\times 10^{-4} & 4.76811\times 10^{-4} & 1.29802\times 10^{-3} & 2.75976\times
   10^{-3} \\
10& 3.65441\times 10^{-5} & 1.59551\times 10^{-4} & 5.04694\times 10^{-4} & 1.21325\times
   10^{-3} \\
11& 1.38971\times 10^{-6} & 5.02895\times 10^{-5} & 1.89440\times 10^{-4} & 5.15815\times
   10^{-4} \\
12& -1.40240\times 10^{-6} & 1.49690\times 10^{-5} & 6.87622\times 10^{-5} & 2.12808\times
   10^{-4} \\
\end{array}
$$
\end{center}
\end{table}

\begin{table}[htbp]
\caption{Fourier coefficients $\widehat{a_{n, k}}^{(50)}\quad (n=20, 22)$}
\label{tab:upptail2022}
\begin{center}
$$
\begin{array}{rrr}
k & n=20 & n = 22  \\
\hline
0& 2.42161\times 10^{-1} & 2.29129\times 10^{-1} \\
1& 2.27078\times 10^{-1} & 2.17130\times 10^{-1} \\
2& 1.88097\times 10^{-1} & 1.85373\times 10^{-1} \\
3& 1.39342\times 10^{-1} & 1.43837\times 10^{-1} \\
4& 9.38340\times 10^{-2} & 1.02652\times 10^{-1} \\
5& 5.83684\times 10^{-2} & 6.82124\times 10^{-2} \\
6& 3.39847\times 10^{-2} & 4.26601\times 10^{-2} \\
7& 1.87134\times 10^{-2} & 2.53303\times 10^{-2} \\
8& 9.82378\times 10^{-3} & 1.43802\times 10^{-2} \\
9& 4.94778\times 10^{-3} & 7.84948\times 10^{-3} \\
10& 2.40299\times 10^{-3} & 4.13869\times 10^{-3} \\
11& 1.13008\times 10^{-3} & 2.11578\times 10^{-3} \\
12& 5.16361\times 10^{-4} & 1.05205\times 10^{-3} \\
\end{array}
$$
\end{center}
\end{table}

\begin{figure}[htbp]
\begin{center}
\includegraphics{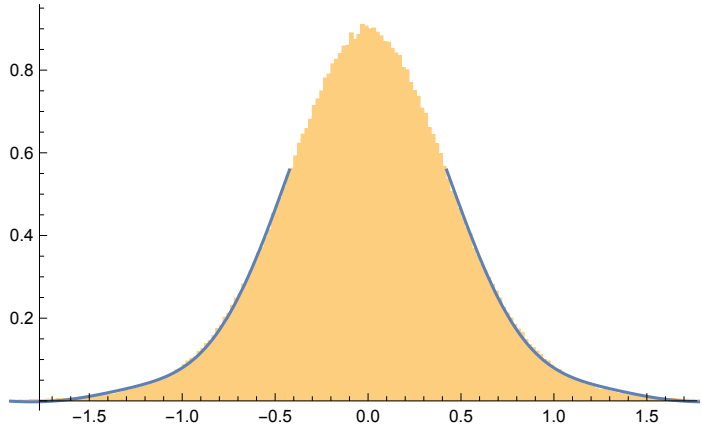}
\end{center}
\caption{Plot of $y=\widehat{f_{20}}^{(12, 50)}(x)$ and 
a histgram of $\surd{b_{1}}$ ($10^6$ replications)}
\label{fig:n20.pdf}
\end{figure}

We now illustrate the case $n=20$ (see also Figure~\ref{fig:n20.pdf}):
\begin{align*}
\widehat{f_{20}}^{(12, 50)}(x) \, = & \,\,
1.21081\times 10^{-1}   
+ \left(2.27078\times 10^{-1}\right) \cos \left(\frac{1}{18} \sqrt{19} \pi 
   x\right)
   \\ &
   +\left(1.88097\times 10^{-1}\right) \cos \left(\frac{1}{9} \sqrt{19} \pi 
   x\right)+\left(1.39342\times 10^{-1}\right) \cos \left(\frac{1}{6} \sqrt{19} \pi 
   x\right)
   \\ &
+\left(9.38340\times 10^{-2}\right) \cos \left(\frac{2}{9} \sqrt{19} \pi 
   x\right)+\left(5.83684\times 10^{-2}\right) \cos \left(\frac{5}{18} \sqrt{19} \pi 
   x\right)
   \\ &
   +\left(3.39847\times 10^{-2}\right) \cos \left(\frac{1}{3} \sqrt{19} \pi 
   x\right)+\left(1.87134\times 10^{-2}\right) \cos \left(\frac{7}{18} \sqrt{19} \pi 
   x\right)
   \\ &
   +\left(9.82378\times 10^{-3}\right) \cos \left(\frac{4}{9} \sqrt{19} \pi 
   x\right)+\left(4.94778\times 10^{-3}\right) \cos \left(\frac{1}{2} \sqrt{19} \pi 
   x\right)
   \\ &
+\left(2.40299\times 10^{-3}\right) \cos \left(\frac{5}{9} \sqrt{19} \pi 
   x\right)+\left(1.13008\times 10^{-3}\right) \cos \left(\frac{11}{18} \sqrt{19} \pi 
   x\right)
   \\ &
      +\left(5.16361\times 10^{-4}\right) \cos \left(\frac{2}{3} \sqrt{19} \pi 
   x\right)
\end{align*}
and 
\noindent
\begin{align*}
\widehat{F_{20}}^{(12, 50)}(x) \, = & \,\,
\frac{1}{2} \left(\frac{\sqrt{19} x}{18}+1\right)+\left(2.98484\times 10^{-1}\right) \sin
   \left(\frac{1}{18} \sqrt{19} \pi  x\right)
   \\ &
   +\left(1.23623\times 10^{-1}\right) \sin
   \left(\frac{1}{9} \sqrt{19} \pi  x\right)+\left(6.10528\times 10^{-2}\right) \sin
   \left(\frac{1}{6} \sqrt{19} \pi  x\right)
     \\ &
   +\left(3.08352\times 10^{-2}\right) \sin
   \left(\frac{2}{9} \sqrt{19} \pi  x\right)+\left(1.53445\times 10^{-2}\right) \sin
   \left(\frac{5}{18} \sqrt{19} \pi  x\right)
   \\ & 
   +\left(7.44524\times 10^{-3}\right) \sin
   \left(\frac{1}{3} \sqrt{19} \pi  x\right)+\left(3.51399\times 10^{-3}\right) \sin
   \left(\frac{7}{18} \sqrt{19} \pi  x\right)
    \\ &
   +\left(1.61412\times 10^{-3}\right) \sin
   \left(\frac{4}{9} \sqrt{19} \pi  x\right)+\left(7.22627\times 10^{-4}\right) \sin
   \left(\frac{1}{2} \sqrt{19} \pi  x\right)
    \\ &
   +\left(3.15863\times 10^{-4}\right) \sin
   \left(\frac{5}{9} \sqrt{19} \pi  x\right)+\left(1.35040\times 10^{-4}\right) \sin
   \left(\frac{11}{18} \sqrt{19} \pi  x\right)
    \\ &
   +\left(5.65611\times 10^{-5}\right) \sin
   \left(\frac{2}{3} \sqrt{19} \pi  x\right).
\end{align*}

Tables~\ref{tab:tailprob} and \ref{tab:upptail} list the upper tail probabilities 
$1 - \widehat{F_{20}}^{(12, 50)}(x)$ 
and percentiles $x_{\alpha}$, respectively, for $\alpha=0.900, 0.950, 0.975, 0.990. 0.995, 0.999$. 
Here, $x_{\alpha}$ was obtained by solving $\widehat{F_{20}}^{(12, 50)}(x_{\alpha}) = \alpha$.
In these tables, we underline the values of the decimal places that differ from those of Mulholland~\cite{mulholland1977null}.

\begin{table}[htbp]
\caption{Upper tail probability $1 - \widehat{F_{n}}^{(12, 50)}(x)$ of $\surd{b_{1}}$}
\label{tab:tailprob}
\begin{center}
$$
\begin{array}{rrrrrr}
 x & n=4 & n=6 & n=8 & n=10 & n=12 
 \\
\hline
 0.1 & 0.41\U{78} & 0.4332 & 0.4311 & 0.4276 & 0.4240 \\
 0.2 & 0.3\U{604} & 0.370\U{5} & 0.3648 & 0.3579 & 0.3511 \\
 0.3 & 0.308\U{3} & 0.312\U{3} & 0.302\U{8} & 0.2932 & 0.2839 \\
 0.4 & 0.26\U{31} & 0.258\U{3} & 0.2464 & 0.2352 & 0.2244 \\
 0.5 & 0.220\U{7} & 0.209\U{3} & 0.1967 & 0.1849 & 0.1736 \\
 0.6 & 0.1818 & 0.165\U{3} & 0.1544 & 0.1427 & 0.1316 \\
 0.7 & 0.1451 & 0.127\U{9} & 0.1194 & 0.1082 & 0.0979 \\
 0.8 & 0.110\U{3} & 0.098\U{6} & 0.0907 & 0.080\U{8} & 0.0717 \\
 0.9 & 0.077\U{6} & 0.075\U{6} & 0.067\U{6} & 0.0594 & 0.0517 \\
 1.0 & 0.04\U{59} & 0.0567 & 0.0496 & 0.0430 & 0.0368 \\
 1.1 & 0.016\U{2} & 0.0414 & 0.0359 & 0.0307 & 0.0258 \\
 1.2 & \text{--} & 0.0291 & 0.0257 & 0.0217 & 0.0179 \\
 1.3 & \text{--} & 0.0193 & 0.0181 & 0.0150 & 0.0123 \\
 1.4 & \text{--} & 0.0118 & 0.0124 & 0.0103 & 0.0083 \\
 1.5 & \text{--} & 0.0063 & 0.0082 & 0.0069 & 0.0056 \\
 1.6 & \text{--} & 0.0026 & 0.0051 & 0.0046 & 0.0037 \\
 1.7 & \text{--} & 0.0006 & 0.0030 & 0.0029 & 0.0024 \\
 1.8 & \text{--} & \text{--} & 0.0016 & 0.0018 & 0.0015 \\
 1.9 & \text{--} & \text{--} & 0.0008 & 0.0010 & 0.0009 \\
 2.0 & \text{--} & \text{--} & 0.0003 & 0.0006 & 0.0005 \\
\end{array}
$$
\end{center}
\end{table}

\begin{table}[htbp]
\caption{Upper percentiles $x_{\alpha}$ of $\surd{b_{1}}$}
\label{tab:upptail}
\begin{center}
$$
\begin{array}{rrrrrrr}
n & \multicolumn{6}{c}{\alpha}
\\
\hline
  & 0.900 & 0.950 & 0.975 & 0.990 & 0.995 & 0.999 
\\
\hline
  4& 0.8305 & 0.9869 & 1.0697 & 1.12\U{10} & 1.13\U{79} & 1.1513 \\
  6& 0.7945 & 1.04\U{12} & 1.2392 & 1.4299 & 1.5306 & 1.6707 \\
  8& 0.7652 & 0.9977 & 1.2080 & 1.4524 & 1.6046 & 1.8668 \\
10& 0.7275 & 0.9539 & 1.1595 & 1.4075 & 1.5785 & 1.9065 \\
12& 0.6931 & 0.9100 & 1.1091 & 1.3532 & 1.5255 & 1.8818 \\
14& 0.6622 & 0.8699 & 1.0614 & 1.2985 & 1.4683 & 1.8315 \\
16& 0.6347 & 0.8337 & 1.0176 & 1.2467 & 1.4121 & 1.7739 \\
18& 0.6102 & 0.8012 & 0.9778 & 1.1988 & 1.3595 & 1.71\U{33} \\
20& 0.5881 & 0.7721 & 0.9419 & 1.15\U{43} & 1.31\U{08} & 1.65\U{20} \\
22& 0.5681 & 0.7459 & 0.90\U{97} & 1.11\U{30} & 1.26\U{37} & 1.60\U{45} \\
\end{array}
$$
\end{center}
\end{table}

\section{Concluding remarks}
We derived approximations to the sampling distributions based on the Fourier cosine series, and constructed 
their coefficients using higher order moments. 
The proposed approximations were applied to 
sums of uniformly distributed random variables and 
sample skewness with small $n$.
In future work, we will derive an approximate distribution function of sample kurtosis under normality, for which the recurrence relations of moments such as those in \eqref{eq:momentsrec}
have been given by Nakagawa et al.~\cite{NHN2016}.

\end{document}